\documentclass[12pt,epsf]{article}
\usepackage{amssymb,amsmath,eucal,amsthm}
\usepackage{graphicx}

\def\bce\{\begin{center}
\def\ece{\end{center}}
\newtheorem{theorem}{Theorem}[section]

\newtheorem{proposition}{Proposition}[section]
\theoremstyle{definition}

\theoremstyle{remark}

\renewcommand{\proofname}{\bf Proof}
\newcommand{\ca}{{\cal A}}
\newcommand{\rr}{{\rm{I \! R}}}
\title{A mathematical approach with delay kernel for the role of the immune response time delay in periodic therapy of the tumors}

\author{\small{G. MIRCEA$^{a}$, M. NEAM\c TU$^{a}$\thanks{Corresponding author}, R.F. HORHAT$^b$, D. OPRI\c S$^{c}$} }
\date{ }

\begin{document}
\maketitle

\begin{tabular}{cccccccc}

\scriptsize{$^{a}$Department of Economic Informatics and Statistics, Faculty of Economics,}\\
\scriptsize{West University of Timi\c soara, str. Pestalozzi, nr. 16A, 300115, Timi\c soara, Romania,}\\
\scriptsize{E-mail:gabriela.mircea@fse.uvt.ro, mihaela.neamtu@fse.uvt.ro,}\\
\scriptsize{$^{b}$ Department of Biophysics and Medical Informatics,}\\
\scriptsize{University of Medicine and Pharmacy, Piata Eftimie Murgu, nr. 3, 300041, Timi\c soara, Romania,}\\
\scriptsize{E-mail: rhorhat@medinfo.umft.ro}\\
\scriptsize{$^{c}$ Department of Applied
Mathematics, Faculty of Mathematics,}\\
\scriptsize{West University of Timi\c soara, Bd. V. Parvan, nr. 4, 300223, Timi\c soara, Romania,}\\
\scriptsize{E-mail: opris@math.uvt.ro}\\

\end{tabular}

\begin{abstract} We consider the model of interaction between the immune system and
tumor cells including a memory function that reflect the influence
of the past states, to simulate the time needed by the latter to
develop a chemical and cell mediated response to the presence of
the tumor. The memory function is called delay kernel. The results
are compared with those from other papers, concluding that the
memory function introduces new instabilities in the system leading
to an uncontrolable growth of the tumor. If the coefficient of the
memory function is used as a bifurcation parameter, it is found
that Hopf bifurcation occurs for kernel. The direction and
stability of the bifurcating periodic solutions are determined.
Some numerical simulations for justifying the theoretical analysis
are also given.

% We find the existence of metastable states induced
%by the treatment, and also of potentially adverse effects of the
%dosage frequency on the stabilization of the tumor.
\end{abstract}

\section{Introduction}
\hspace{0.5cm} As everyone knows, cancer is one of the most
fearsome illness. It was declared the disease of  20th century.
Many efforts were made to cure it, but to do this, first of all it
is needed to understand its physiopathological me\-cha\-nisms. It
was discovered that the human body is not completely helpless
against this disease and it fights against cancer using its best
and powerful weapon, namely immune system.

In what follows, we will not make a general overview of the immune
system, but we will mention briefly some of its components and
aspects of the dynamics which appear in our model. The cell who
performs directly the tumor elimination is T-lymphocyte, which is
activated by b-lymphocytes through the cytokines [3,11]. It is not
of less importance to mention the immunodepression, a phenomenon
that appears in the tumor region, when the tumor increases its
size, and leads to the deactivation of the lymphocytes.

In the effort of modeling this process an important role plays the
time delay. It is obvious for everyone that the biological process
do not take place instantaneously and an amount of time is needed
for, in our case, the interaction between immune system and the
tumor [2,4]. During the years, some models, concerning tumor
dynamics have been develop [6,7,12] and some of them includes time
delay [1,10,13].

In what follows we propose a model of the interaction tumor-immune
system using delay kernel.

Let $x(t)$ and $y(t)$ denote respectively the number of malignant
and lymphocyte cells, for $t\in \mathbb{R}$. The rate of malignant
cells $(\dot{x}(t))$ is given by [12]:
\begin{equation}
\dot{x}(t)=a_{1}x(t)-a_{2}x(t)y(t).
\end{equation}
We assume that the growth rate is proportional to $x(t)$ and the
decrease rate is proportional to the frequency of interaction with
lymphocytes. The coefficients are $a_{1}$ and  $a_{2}$,
respectively,
where $a_{1}$ is tissue dependent.\\
On the other hand, the growth rate of lymphocytes $\dot{y}(t)$ is
described by [12]:
\begin{equation}
\dot{y}(t)=b_{1}x(t)y(t)-b_{2}x(t)-b_{3}y(t)+b_{4}.
\end{equation}
It is proportional to the interaction with malignant cells and
also to the flux per unit time of lymphocytes to the place of
interaction. These effects are represented by the first and fourth
terms in the right-hand side of equation (2). The mortality of the
lymphocytes is proportional with y(t) (natural death) and also
x(t), which express the immunodepression phenomenon. The term
$b_{1}x(t)y(t)$ is important for this study. It means the
interaction between the two populations, $x(t)$ and $y(t)$ with a
frequency $b_{1}$ of recognition of malignant cells by the immune
system. We consider the effect of influence of the past for this
chemical signal mediated interaction which introduces the memory
functions $\rho_{1}$ and $\rho_{2}$, which are nonnegative bounded
functions defined on $[0,\infty)$ and
\begin{eqnarray*} \\
\int\limits_{0}^{\infty}k_i(s)ds=1,\quad
\int\limits_{0}^{\infty}sk_i(s)ds< \infty, i=1,2.
\end{eqnarray*}
The evolution equations (1), (2) become now
\begin{equation}
\begin{split}
&\dot{x}(t)=a_{1}x(t)-a_{2}x(t)y(t)\\
&\dot{y}(t)=b_{1}(\int\limits_{0}^{\infty}k_{1}(s)x(t-s)ds)
(\int\limits_{0}^{\infty}k_{2}(s)y(t-s)ds)-b_{2}x(t)-b_{3}y(t)+b_{4}.
\end{split}
\end{equation}
The memory functions are called delay kernels. The delay becomes a
discrete one when the delay kernel is a delta function at a
certain time. Usually, we employ the following form
\begin{equation*}
k_{i}(s)=\frac{1}{p!}q_{i}^{p+1}s^{p}e^{-q_{i}s},\quad i=1,2,
\end{equation*}for the memory function.
When $p=0$ and $p=1$, the memory functions are called "weak"
and "strong" kernel, respectively. \\
For $k_{i}(s)=\delta(s-\tau_{i}), i=1,2, \tau_{1}\geq
0,\tau_{2}\geq 0$ equation (4) is given by
\begin{equation}
\begin{split}
&\dot{x}(t)=a_{1}x(t)-a_{2}x(t)y(t)\\
&\dot{y}(t)=b_{1}x(t-\tau_{1})y(t-\tau_{2})-b_{2}x(t)-b_{3}y(t)+b_{4}.
\end{split}
\end{equation}
The model (4) with $\tau_{1}=\tau_{2}=\tau$, is the model from [4]
which has been studied using only numerical simulations.\\
In this paper, we analyze the model (4) with the following initial
values
\begin{equation*}
x_{1}(\theta)=\varphi_{1}(\theta),\quad
x_{2}(\theta)=\varphi_{2}(\theta), \quad \theta\in(-\infty,0]
\end{equation*}
and $\varphi_{1},\varphi_{1} $ as differentiable functions.\\
The paper is organized as follows. In section 2, we discuss the
local sta\-bi\-li\-ty for the equilibrium states of system (4),
for different forms of the delay kernels. We investigate the
existence of the Hopf bifurcation with respect of the parameters
of the delay kernels. In section 3, the direction of the Hopf
bifurcation is analyzed by normal form theory and the center
manifold theorem. Numerical simulations in order to justify the
theoretical results are illustrated in section 4. Finally, some
conclusions are made.

\section{Local stability and existence of the Hopf bifurcation}
We consider model (4) with parameters $a_{1}, a_{2}, b_{1}, b_{2},
b_{3}, b_{4}$ assumed positives numbers and $\frac{b_{2}}{b_{1}}<
\frac{b_{4}}{b_{3}}< \frac{a_{1}}{a_{2}}$. The equilibrium states
of system (4) are the points $L_{0}=(x_{0},y_{0})$ and
$L_{1}=(0,\frac{b_{4}}{b_{3}})$, where
\begin{equation*}
x_{0}=\frac{b_{3}a_{1}-b_{4}a_{2}}{a_{1}b_{1}-a_{2}b_{2}}, \quad
y_{0}=\frac{a_{1}}{a_{2}}.
\end{equation*}
We analyzed the  local stability in the equilibrium state $L_{0}$.
We consider the following translation
\begin{equation}
x_{1}(t)=x(t)-x_{0}, \quad x_{2}(t)=y(t)-y_{0}.
\end{equation}
With respect to (5), the system (4) can be expressed as
\begin{equation}
\begin{split}
& \dot{x}_{1}(t)=-a_{2}x_{0}x_{2}(t)-a_{2}x_{1}(t)x_{2}(t)\\
&\dot{x}_{2}(t)=-b_{2}x_{1}(t)-b_{3}x_{2}(t)+b_{1}x_{0}
\int\limits_{0}^{\infty}k_{2}(s)x_{2}(t-s)ds+b_{1}y_{0}
\int\limits_{0}^{\infty}k_{1}(s)x_{1}(t-s)ds\\&+
b_{1}(\int\limits_{0}^{\infty}k_{1}(s)x_{1}(t-s)ds)
(\int\limits_{0}^{\infty}k_{2}(s)x_{2}(t-s)ds).
\end{split}
\end{equation}
The system (6) has  0=(0,0) as equilibrium state.\\
To investigate the local stability of equilibrium state of the
system (8), we linearize system (6). The linearized system of (6)
is
\begin{equation}\dot{U}(t)=AU(t)+B_{1}U_{1}(t)+B_{2}U_{2}(t),
\end{equation}
where
\begin{equation}
\begin{split}
A=\left(\begin{array}{ccc} 0&-a_{2}x_{0}\\-b_{2}&-b_{3}
\end{array} \right), B_{1}=\left(\begin{array}{ccc}0&0\\
b_{1}y_{0}& 0 \end{array} \right),B_{2}=\left(\begin{array}{ccc}
0&0\\0 & b_{1}x_{0}\end{array} \right)
\end{split}
\end{equation}
with
\begin{equation*}U(t)=(u_{1}(t),u_{2}(t))^{T},\,U_{i}(t)=
(\int\limits_{0}^{\infty}k_{i}(s)u_{1}(t-s)ds,
\int\limits_{0}^{\infty}k_{i}(s)u_{2}(t-s)ds)^{T}, \, i=1,2.
\end{equation*}
The characteristic equation corresponding to system (7) is
$\Delta(\lambda)$=0, where
\begin{equation}
\Delta(\lambda)=det(\lambda
I-A-(\int\limits_{0}^{\infty}k_{1}(s)e^{-\lambda s}ds)B_{1}-
(\int\limits_{0}^{\infty}k_{2}(s)e^{-\lambda s}ds)B_{2}).
\end{equation}
From (8) and (9), we have:
\begin{equation}
\Delta(\lambda)=\lambda^{2}+b_{3}\lambda-a_{2}b_{2}x_{0}+a_{1}b_{1}x_{0}
\int\limits_{0}^{\infty}k_{1}(s)e^{-\lambda s}ds-\lambda
b_{1}x_{0}\int\limits_{0}^{\infty}k_{2}(s)e^{-\lambda s}ds.
\end{equation}
The equilibrium state $L_{0}$ is locally asymptotically stable if
and only if the eigenvalues of $\Delta(\lambda)=0$ have negative
real
parts.\\
Because of the presence of two delay kernels $k_{1}$ and $k_{2}$
in the equation $\Delta(\lambda)=0$, the analysis of the sign of
real parts of eigenvalues is complicated and a direct
approach cannot be considered.\\
We analyze the eigenvalues for the equation $\Delta(\lambda)=0$ if
the delay kernels $k_{1}$ and $k_{2}$ are delta functions
or $k_1$ is delta function and $k_2$ is weak function.\\
Using results from [3], we obtain:
\begin{proposition}\label{prop1}
If
\begin{equation}
k_{1}(s)=\delta(s-\tau_{1}),\, k_{2}(s)=\delta(s-\tau_{2}),\,
\tau_{1}\geq 0,\tau_{2}\geq 0
\end{equation}then
\begin{itemize}
\item[(i)] function $(10)$ is given by
\begin{equation}\Delta(\lambda,\tau_{1},\tau_{2})=
\lambda^{2}+b_{3}\lambda-a_{2}b_{2}x_{0}+a_{1}b_{1}x_{0}e^{-\lambda
\tau_{1}}- \lambda b_{1}x_{0}e^{-\lambda\tau_{2}};
\end{equation}
\item[(ii)]if $\tau_{1}=0, \tau_{2}=0$ then the equilibrium state
$L_{0}$ of system $(4)$ is locally asymptotic stable;
 \item[(iii)]if
 \begin{equation*}0\leq
 \tau_{1}+\tau_{2}<\frac{b_{3}+b_{1}x_{0}}{a_{1}b_{1}x_{0}},
 \end{equation*}
 then the equilibrium state $L_{0}$ of the system $(4)$ is
 asymptotically stable.
\end{itemize}
\end{proposition}
Next, we study the existence of Hopf bifurcation of system (4)
with $k_{1}$ and $k_{2}$ given by (11), by choosing one of the
delays as a bifurcation parameter, e.g. take $\tau_{1}$ as the
bifurcation parameter. First, we would like to know when
$\Delta(\lambda,\tau_{1},\tau_{2})=0$, where
$\Delta(\lambda,\tau_{1},\tau_{2})$ given by (12) has purely
imaginary roots $\lambda=\pm i\omega_{0}(\omega_{0}> 0)$ at
$\tau_{1}=\tau_{10}$. Note that
\begin{equation}
\begin{split}
&\omega_{0}^{2}+a_{2}b_{2}x_{0}-a_{1}b_{1}x_{0}\cos(\tau_{10}\omega_{0})
+b_{1}x_{0}\omega_{0}\sin(\tau_{2}\omega_{0})=0\\
 & b_{3}\omega_{0}-a_{1}b_{1}x_{0}\sin(\tau_{10}\omega_{0})-
 b_{1}x_{0}\omega_{0}\cos(\tau_{2}\omega_{0})=0,
\end{split}
\end{equation}
which implies that
\begin{equation*}
\sin((\tau_{10}-\tau_{2})\omega_{0})=g(\omega_{0}),
\end{equation*}
where
\begin{equation}g(\omega)=\frac{\omega^{4}-(b_{1}^{2}x_{0}^{2}-b_{3}^{2}
-2a_{2}b_{2}x_{0})\omega^{2}+(a_{2}^{2}b_{2}^{2}-a_{1}^{2}b_{1}^{2})x_{0}^{2}}{
2a_{1}b_{1}^{2}x_{0}^{2}\omega}.
\end{equation}
From (14),  $g'(\omega)>0$. So $g(\omega)$ is strictly
monotonically increasing on $[0,\infty)$, with
$\lim\limits_{\omega\rightarrow 0}g(\omega)=-\infty $ and
$\lim\limits_{\omega\rightarrow \infty}g(\omega)=\infty $.
Clearly, if $\tau_{10}>\tau_2$ then $g(\omega)$ intersects
$\sin((\tau_{10}-\tau_{2})\omega)$ only in a point. Hence
$\lambda=i\omega_{0}$ is a simple root of equation
$\Delta(\lambda,\tau_{1},\tau_{2})=0$. Differentiating
$\Delta(\lambda,\tau_{1},\tau_{2})$ implicitly with respect to
$\tau_{1}$, we obtain
\begin{equation*}
\begin{split}
& Re[\frac{d\lambda}{d\tau_{1}}]_{\lambda=i\omega_{0},
\tau_{1}=\tau_{10}}=\\
& =-\frac{a_{1}^{2}b_{1}^{2}x_{0}^{2}
\!+\!a_{1}b_{1}^{2}x_{0}^{2}\tau_{2}\omega_{0}^{2}
\cos((\tau_{10}\!-\!\tau_{2})\omega_{0})\!+\!a_{1}b_{1}x_{0}^{2}
(\omega_{0}^{2}\!+\!a_{2}b_{2}x_{0})\cos(\tau_{10}\omega_{0})}{l_{1}^{2}+l_{2}^{2}}
\end{split}
\end{equation*}
where
\begin{equation*}
\begin{split}
& l_{1}=b_{3}-a_{1}b_{1}x_{0}\tau_{10}\cos(\tau_{10}\omega_{0})
-b_{1}x_{0}\cos(\tau_{2}\omega_{0})+b_{1}x_{0}\tau_{2}\omega_{0}
\sin(\tau_{2}\omega_{0})\\
&l_{2}=2\omega_{0}-a_{1}b_{1}x_{0}\tau_{10}\sin(\tau_{10}\omega_{0})
+b_{1}x_{0}\sin(\tau_{2}\omega_{0})+b_{1}x_{0}\tau_{2}\omega_{0}
\cos(\tau_{2}\omega_{0}).
\end{split}
\end{equation*}
From the above analysis and the standard Hopf bifurcation theory,
we have the following result:
\begin{proposition}\label{prop2}
If $k_{1}(s)=\delta(s-\tau_{1}),\, k_{2}(s)=\delta(s-\tau_{2})$
and there is $\tau_{1}=\tau_{10}$ for given $\tau_2>0$,
$\tau_{10}>\tau_2$ so that equations $(13)$ hold and
\begin{equation*}
Re(\frac{d\lambda}{d\tau_{1}})_{\lambda=i\omega_{0},
\tau_{1}=\tau_{10}}\neq 0,
\end{equation*}
then a Hopf bifurcation occurs at $L_{0}$ as $\tau_{1}$ passes
through $\tau_{10}$.\\
For given $\tau_{2}> 0$, a solution for $(13)$ is the pair
$(\tau_{10},\omega_{10})$, where
\begin{equation*}
\tau_{10}=\frac{k\pi}{\omega_{10}}+\tau_{2}, \, k=1,2,\dots
\end{equation*}
and $\omega_{10}$ is a positive root of the equation
\begin{equation*}
x^{4}-(b_{1}^{2}x_{0}^{2}-b_{3}^{2}-2a_{2}b_{2}x_{0})x^{2}
+(a_{2}^{2}b_{2}^{2}-a_{1}^{2}b_{1}^{2})x_{0}^{2}=0.
\end{equation*}
\end{proposition}

\begin{proposition}\label{prop3}If
\begin{equation}
k_{1}(s)=\delta(s-\tau_{1}),\, k_{2}(s)=q_{2}e^{-q_{2}s},
\tau_{1}\geq 0,q_{2}>0
\end{equation}then
\begin{itemize}
\item[(i)] function $(10)$ is given by
\begin{equation}\Delta(\lambda,\tau_{1},q_{2})=\lambda^{3}+p_{2}\lambda^{2}
+p_{1}\lambda+p_{0}+(r_{1}\lambda+r_{0})e^{-\lambda\tau_{1}},
\end{equation}where
\begin{equation}
\begin{split}
&p_{2}= q_{2}+b_{3}, \,
p_{1}=q_{2}b_{3}-a_{2}b_{2}x_{0}-b_{1}x_{0}q_{2}\\
& p_{0}=-q_{2}a_{2}b_{2}x_{0}, \, r_{1}=a_{1}b_{1}x_{0},\,
r_{0}=a_{1}b_{1}x_{0}q_{2};
\end{split}
\end{equation}
\item[(ii)]if $\tau_{1}=0$ and
\begin{equation*}
4(a_{1}b_{1}-a_{2}b_{2})^{2}<a_{2}b_{3}(b_{1}b_{4}-a_{2}b_{3})
\end{equation*}
then for $q_{2}\in (0,q_{21})\cup (q_{22},\infty)$ the equilibrium
state $L_{0}$ is locally asymptotic stable, where $q_{21},q_{22}$
are the solutions of the equation
\begin{equation*}
(b_{3}-b_{1}x_{0})x^{2}+b_{3}(b_{3}-b_{1}x_{0})x+b_{3}(a_{1}b_{1}-a_{2}b_{2})x_{0}=0.
\end{equation*}
\end{itemize}
\end{proposition}
Next, we study the existence of Hopf bifurcation for system $(3)$
with $\rho_{1}$ and $\rho_{2}$ given by (15), by choosing the
delay $\tau_{1}$ as the bifurcation parameter. First, we  would
like to know when $\Delta(\lambda,\tau_{1},q_{2})=0$, where
$\Delta(\lambda,\tau_{1},q_{2})$ is given by (16), has purely
imaginary roots $\lambda=\pm i\omega_{01}(\omega_{01}> 0)$ at
$\tau_{1}=\tau_{11}$. Note that
\begin{equation}
\begin{split}
& p_{0}-p_{2}\omega_{01}^{2}+r_{0}\cos(\omega_{01}\tau_{11})+
r_{1}\sin(\omega_{01}\tau_{11})=0\\
&-\omega_{01}^{3}+p_{1}\omega_{01}+r_{1}\omega_{01}\cos(\omega_{01}\tau_{11})-
r_{0}\sin(\omega_{01}\tau_{11})=0,
\end{split}
\end{equation}
which implies that
\begin{equation}
\omega_{01}^{6}+(p_{2}^{2}-2p_{1})\omega_{01}^{4}+
(p_{1}^{2}-2p_{0}p_{2}+r_{1}^{2})\omega_{01}^{2}+p_{0}^{2}-r_{0}^{2}=0.
\end{equation}
From (17), $p_{0}^{2}<r_{0}^{2}$ and from (19) $\lambda=
i\omega_{01}$ is a simple root of the equation
$\Delta(\lambda,\tau_{1},q_{2})=0$. From (18) we obtain:
\begin{equation*}
\tau_{11}=\frac{1}{\omega_{01}}\arctan\frac{r_{1}\omega_{01}(p_{2}\omega_{01}^{2}-p_{0})+r_{0}(p_{1}
\omega_{01}-\omega_{01}^{3})}{p_{1}\omega_{01}(\omega_{01}^{3}-p_{1}\omega_{01}
)+r_{0}(p_{0}-p_{2}\omega_{01}^{2})}.
\end{equation*}
Differentiating $\Delta(\lambda,\tau_{1},q_{2})=0$ implicitly with
respect to $\tau_{1}$, we obtain
\begin{equation*}
Re(\frac{d\lambda}{d\tau_{1}})_{\lambda=i\omega_{01},\tau_{1}=\tau_{11}}=
\frac{\omega_{01}(r_{1}\omega_{01}l_{1}-l_{2}r_{0})}{m_{1}^{2}+m_{2}^{2}},
\end{equation*}
where
\begin{equation*}
\begin{split}
&m_{1}=(p_{1}-3\omega_{01}^{2})\cos(\omega_{01}\tau_{11})-2p_{2}\omega_{01}
\sin(\omega_{01}\tau_{11})+r_{1}-r_{0}\tau_{11}\\
&m_{2}=2p_{2}\omega_{01}\cos(\omega_{01}\tau_{11})+(p_{1}-3\omega_{01}^{2})
\sin(\omega_{01}\tau_{11})-r_{1}\tau_{11}-r_{1}\tau_{11}\omega_{01}.
\end{split}
\end{equation*}
\begin{proposition}\label{prop4} If $k_{1}(s)=\delta(s-\tau_{1}),\, k_{2}(s)=q_{2}e^{-q_{2}s}$ and
$\tau_{1}=\tau_{11}$ then
\begin{equation*}
Re(\frac{d\lambda}{d\tau_{1}})_{\lambda=i\omega_{01},\tau_{1}=\tau_{11}}\neq
0
\end{equation*}
and a Hopf bifurcation occurs at $L_{0}$ as $\tau_{1}$ passes
through $\tau_{11}$.
\end{proposition}

\section{Direction and stability of the Hopf bifurcation for
$k_{1}(s)=\delta(s-\tau_{1})$, $k_{2}(s)=\delta(s-\tau_{2})$}

In what follows, we will study the direction and stability in two
cases: in the first case the both kernels are delta function and
in the second case the kernel $k_1$ is delta function and the
kernel $k_2$ is weak function.

\vspace{0.2cm} {\bf 3.1. The case $k_1(s)=\delta(s-\tau_1)$,
$k_2(s)=\delta(s-\tau_2)$, $\tau_1\geq0$, $\tau_2\geq0$.}

In Proposition \ref{prop1} and \ref{prop2}, we obtained some
conditions which guarantee that system (4) undergoes Hopf
bifurcation at $\tau=\tau_{10}$. In this section, we study the
direction, the stability and the period of bifurcating periodic
solutions. The method that we used is based on the normal form
theory
and the center manifold theorem introduced by [2].\\
From the previous section, we know that if $\tau=\tau_{10}$, then
all the roots of $\Delta(\lambda,\tau_{10},\tau_{2})=0$, other
than $\pm i\omega_{0}$ have negative real parts and any root of
the form $\lambda(\tau_{1})=\alpha(\tau_{1})+i\omega(\tau_{1})$
satisfies $\alpha(\tau_{10})=0$, $\omega(\tau_{10})=\omega_{0}$
and $\frac{d\alpha(\tau_{10})}{d\tau_{1}}\neq 0$. For notational
convenience let $\tau_{1}=\tau_{10}+\mu, \mu\in\mathbb{R}$. Then
$\mu =0$ is the Hopf bifurcation value for (4). Without loss of
generality, assume that $\tau_{10}>\tau_{2}$ and define the space
of $C^{1}$ functions as
$C^{1}=C^{1}([-\tau_{10},0],\mathbb{C}^{2})$.\\
Suppose that for given $a_{1}, a_{2},b_{1}, b_{2}, b_{3},
b_{4},\tau_{2}$, there is a $\tau_{10}>0$ at which (4) exhibits a
Hopf bifurcation. In $\tau_{1}=\tau_{10}+\mu, \mu\in\mathbb{R}$,
we regard $\mu$ as the bifurcation parameter. For $\phi\in C^{1}$,
we define

\begin{equation*}
\ca(\mu)\phi=\left\{\begin{array}{rl}
\frac{d\phi(\theta)}{d\theta},&\theta\in [-\tau_{10},0)\\
\int\limits_{-\tau_{10}}^{0}d\eta(t,\mu)\phi(t),&
\theta=0\end{array} \right.
\end{equation*}and

\begin{equation*}
R(\mu)\phi=\left\{\begin{array}{rl} (0, 0)^T, &\theta\in [-\tau_{10},0)\\
(-a_2\phi_1(0)\phi_2(0),
b_{1}\phi_{1}(-\tau_{10})\phi_{2}(-\tau_{2}))^T,& \theta=0
\end{array}\right.,
\end{equation*}
where

\begin{equation*}
\eta(\theta,\mu)=\left\{\begin{array}{rl}
A,&\theta=0\\
B_{1}\delta(\theta+\tau_{2}),& \theta\in
[-\tau_{2},0)\\
-B_{2}\delta(\theta+\tau_{10}), &\theta\in [-\tau_{10},-\tau_{2})
\end{array} \right.
\end{equation*}and
$A, B_{1}, B_{2}$ are given by (8).\\ Then, we can rewrite (4) in
the following vector form

\begin{equation*}\dot{U}_{t}=\ca(\mu)U_{t}+RU_{t},
\end{equation*}where
\begin{equation*}U=(u_{1},u_{2})^{T},\quad
U_{t}=U(t+\theta),\quad \theta\in [-\tau_{10},0].
\end{equation*}
For $\psi \in C^{1}([0,\tau_{10}],\mathbb{C}^{2})$, the adjoint
operator $\ca^{\ast}$ of $\ca$ is defined as

\begin{equation*}\ca^{\ast}\psi(s)=\left\{\begin{array}{rl}
-\frac{d\psi(s)}{ds},& s\in (0,\tau_{10}] \\
\int\limits_{-\tau_{10}}^{0}d\eta^{T}(t,0)\psi(-t), & s=0
\end{array}\right.
\end{equation*}

For $\phi\in C([-\tau_{10},0],\mathbb{C}^{2})$ and $\psi \in
C([0,\tau_{10}],\mathbb{C}^{2})$ we define the bilinear form

\begin{equation}<\phi,\psi>=\overline{\psi}^{T}(0)\phi(0)-
\int\limits_{-\tau_{10}}^{0}\int\limits_{\xi=0}^{\theta}
\overline{\psi}^{T}(\xi-\theta) d\eta(\theta)\phi(\xi)d \xi,
\end{equation}
where $\eta(\theta)=\eta(\theta,0)$.

\begin{proposition}\label{prop5}
\begin{itemize}
\item[(i)]The eigenvector  of $\ca(0)$ corresponding to eigenvalue
$\lambda_{1}=i\omega_{0}$ is given by
\begin{equation*}
h(\theta)=(v_{1},v_{2})^{T}e^{\lambda_{1}\theta}, \quad
\theta\in[-\tau_{10},0],
\end{equation*} where

\begin{equation*}v_{1}=1,\quad
v_{2}=\frac{a_{2}b_{2}-a_{1}b_{1}e^{\lambda_{2}\tau_{10}}}
{a_2(b_{3}+\lambda_{1}-b_{1}x_{0}e^{\lambda_{2}\tau_{2}})}
\end{equation*}
and $\lambda_{2}=\overline{\lambda_{1}}$; \item[(ii)]The
eigenvector  of $\ca^{\ast}$ corresponding to  eigenvalue
$\lambda_{2}$ is
\begin{equation*}h^{\ast}(s)=(w_{1},w_{2})^{T}e^{\lambda_{1}s},
\quad s\in [0,\infty)
\end{equation*}
where
\begin{equation*}
\begin{split}
& w_{1}=\frac{f_1}{\eta},\quad w_{2}=\frac{1}{\eta}, \quad f_1=\frac{a_{2}b_{2}-a_{1}b_{1}e^{\lambda_{1}\tau_{10}}}{a_2\lambda_{1}}\\
& \eta=(f_1+b_{1}y_{0}\tau_{10}e^{\lambda_{1}\tau_{10}})+
\overline{v_{2}}(1+\tau_{2}b_{2}x_{0}e^{\lambda_{1}\tau_{2}});
\end{split}
\end{equation*}
\item[(iii)]With respect of (20), we have:
\begin{equation*}<h^{\ast},h>=1,\quad
<h^{\ast},\overline{h}>=<\overline h^{\ast},h>=0, \quad <\overline
h^{\ast},\overline{h}>=1.
\end{equation*}
\end{itemize}
\end{proposition}

Using the approach of [2], [9] we next compute the coordinates of
the center manifold $\Omega_{0}$ at $\mu=0$. Let
$X_{t}=X(t+\theta), \theta\in [-\tau_{10},0)$ be the solution of
system (3) when $\mu =0$.\\
Define

\begin{equation*}
z(t)=<h^{\ast},X_{t}>\quad  w(t,\theta)=X_{t}- 2 Re
(z(t)h(\theta)).
\end{equation*}

On the center  manifold $\Omega_{0}$, we have
\begin{equation*}w(t,\theta)=w(z(t),\overline{z}(t),\theta)=
w_{20}(\theta)\frac{z^{2}}{2}+w_{11}(\theta)
z\overline{z}+w_{02}(\theta)\frac{\overline{z}^{2}}{2}+\cdots
\end{equation*}
where $z$ and $\overline{z}$ are the local coordinates of the
center manifold $\Omega_{0}$ in the direction of $h$ and
$h^{\ast}$, respectively.\\
For the solution $u_{t}\in\Omega_{0}$, we have:

\begin{equation*}\dot{z}(t)=\lambda_{1}z(t)+g(z(t),\overline{z(t)}),
\end{equation*}where
\begin{equation*}g(z,\overline{z})=g_{20}\frac{z^{2}}{2}+g_{11}z\overline{z}+g_{21}
\frac{\overline{z}^{2}}{2}+g_{21}\frac{z^{2}\overline{z}}{2}
\end{equation*}

\begin{proposition}For the system (4), the coefficients $g_{20}, g_{11},
 g_{02}, g_{21}$ and the functions $w_{20}(\theta), w_{11}(\theta), w_{02}(\theta)$
are given by
\begin{equation}
\begin{split}
& g_{20}=\overline w_1f_{120}+\overline w_2f_{220},\quad
g_{11}=\overline w_1f_{111}+\overline w_2f_{211},\\
& g_{02}=\overline w_1f_{102}+\overline w_2f_{202},\quad
g_{21}=\overline w_1f_{121}+\overline w_2f_{221}.
\end{split}
\end{equation}
where
\begin{equation*}
\begin{split}
&  f_{120}=-2a_2v_{1}v_{2}, f_{111}=-2a_2Re(v_1\overline{v_{2}}),
f_{102}=\overline
f_{120},\\
& f_{220}=2b_{1}v_1v_{2}e^{\lambda_{2}(\tau_{10}+\tau_{2})},
f_{211}=2b_{1}Re(v_1\overline{v_{2}}e^{\lambda_{1}\tau_{2}+\lambda_{2}\tau_{10}}),
f_{202}=\overline f_{220}\\
& f_{121}=-a_{2}(2v_{1}w_{211}(0)+\overline{v_{1}}w_{220}(0)+2v_{2}w_{111}(0)+\overline{v_{1}}w_{120}(0))\\
&
f_{221}=b_{1}(2v_{1}e^{\lambda_{2}\tau_{10}}w_{211}(-\tau_2)+\overline{v_{1}}e^{\lambda_{1}\tau_{10}}w_{220}(-\tau_2)+
2v_{2}e^{\lambda_{2}\tau_{2}}w_{111}(-\tau_{10})+\\
&+\overline{v_{2}}e^{\lambda_{1}\tau_{2}}w_{120}(-\tau_{10}))
\end{split}
\end{equation*}
and
\begin{equation*}
w_{20}(\theta)=(w_{120}(\theta), w_{220}(\theta))^T, \quad
w_{11}(\theta)=(w_{111}(\theta), w_{211}(\theta))^T,
\end{equation*}

\begin{equation*}
\begin{split}
&w_{20}(\theta)=-\frac{g_{20}}{\lambda_{1}}h(0)e^{\lambda_{1}\theta}
-\frac{\overline{g_{20}}}{3\lambda_{1}}\overline{h(0)}e^{\lambda_{2}\theta}+
E_{1}e^{2\lambda_{1}\theta}\\
&w_{11}(\theta)=\frac{g_{11}}{\lambda_{1}}h(0)e^{\lambda_{1}\theta}-
\frac{\overline{g_{11}}}{\lambda_{1}}\overline{h(0)}e^{\lambda_{2}\theta}+E_{2},
\end{split}
\end{equation*}
and
\begin{equation*}
E_{1}=(E_{11},E_{12})^{T}, \quad E_{2}=(E_{21},E_{22})^{T},
\end{equation*}
where
\begin{equation*}
\begin{split}
&E_{11}=\frac{(2\lambda_{1}+b_{3}-
b_{1}x_{0}e^{\lambda_{1}\tau_{2}})f_{120}-a_2x_0f_{220}}{2\lambda_{1}(-2\lambda_{1}-b_{3}+
b_{1}x_{0}e^{2\lambda_{1}\tau_{2}})+a_2x_0(-b_2+b_1y_0e^{2\lambda_{1}\tau_{10}})}\\
&E_{12}=\frac{2\lambda_1E_{11}+f_{120}}{a_2x_0},
E_{21}=\frac{(b_1x_0-b_3)E_{22}+f_{211}}{b_1y_0-b_2},
E_{22}=-\frac{f_{111}}{a_2x_0}.
\end{split}
\end{equation*}
\end{proposition}

Based on the above analysis and calculation, we can see that each
$g_{ij}$ in (21) is determined by the parameters and delays in
system (3). Thus, we can explicitly compute the following
quantities:
\begin{equation}
\begin{split}&C_{10}(0)=\frac{i}{2\omega_{0}}(g_{20}g_{11}-2|g_{11}|^{2}-\frac{1}{3}
|g_{02}|^{2})+\frac{g_{21}}{2}\\
&\mu_{20}=-\frac{Re\, C_{10}(0)}{Re \,\lambda'(0)}\\
&T_{20}=-\frac{Im\,C_{10}(0)+\mu_{20}Im\,\lambda'(0)
}{\omega_{0}}\\
&\beta_{20}=2Re\, (C_{10}(0)).
\end{split}
\end{equation}

In summary, this leads to the following result:
\begin{theorem}In formulas (22), $\mu_{20}$ determines the directions of the Hopf
bifurcations: if $\mu_{20}>0(<0)$ the Hopf bifurcation is
supercritical (subcritical) and the bifurcating periodic solutions
exist for $\tau_{1}>\tau_{10}(<\tau_{10})$; $\beta_{20}$
determines the stability of the bifurcation periodic solutions:
the solutions are orbitally stable (unstable) if $\beta_{20}< 0 (>
0)$ and $T_{20}$ determines the periodic solutions: the period
increases (decreases)
 if $T_{20}>0 (<0)$. \\
In $(22)$ $Re (\lambda'(0))$ and $Im(\lambda'(0))$ are given by
\begin{equation*}
\begin{split}
&Re(\lambda'(0))=Re(\frac{d\lambda}{d\tau_{1}})_{\lambda=i\omega_{0},
\tau=\tau_{10}}\\
&
Im(\lambda'(0))=Im(\frac{d\lambda}{d\tau_{1}})_{\lambda=i\omega_{0},
\tau=\tau_{10}}
\end{split}
\end{equation*}
where
\begin{equation*}
\frac{d\lambda}{d\tau_{1}}=\frac{a_{1}b_{1}x_{0}\lambda
e^{-\lambda\tau_{1}}}{b_{3}+2\lambda-a_{1}b_{1}x_{0}\tau_{1}
e^{-\lambda\tau_{1}}-b_{1}x_{0}(1-\lambda\tau_{2})e^{-\lambda\tau_{2}}}.
\end{equation*}
\end{theorem}

\medskip

{\bf 3.2. The case $k_1(s)=\delta(s-\tau_1)$,
$k_2(s)=q_2e^{-q_2s}$, $\tau_1\geq0$, $q_2>0$.}

\hspace{0.6cm} For  $k_1(s)=\delta(s-\tau_1)$,
$k_2(s)=q_2e^{-q_2s}$, $\tau_1\geq0$, $q_2>0$, system (6) is given
by:

\begin{equation}\begin{array}{l} \vspace{0.1cm}
\dot x{}_1(t)=-a_2x_0x_2(t)-a_2x_1(t)x_2(t),\\
\vspace{0.1cm}
\dot x{}_2(t)=-b_2x_1(t)\!-\!b_3x_2(t)\!+\!b_{1}x_0x_{3}(t)\!+\!b_{1}y_0x_2(t\!-\!\tau_1)\!+\!b_1x_3(t)x_2(t\!-\!\tau_1),\\
\dot x{}_3(t)=q_2(x_2(t)-x_3(t)) .\end{array}\end{equation}

We linearize system (23) and obtain:
$$\dot V(t)=A_{1}V(t)+C_{1}V(t-\tau_1),$$
where
$$A_{1}\!\!=\!\!\left(\!\!\!\!\begin{array}{ccc}
\vspace{0.2cm}
0 & -a_2x_0 & 0 \\
\vspace{0.2cm}
-b_2 & -b_3 & b_{1}x_{0}\\
0 & q_2 &  -q_2\\
\end{array}\right)$$
$$C_{1}=\left(\begin{array}{ccc} \vspace{0.2cm}
0 & 0 & 0 \\
\vspace{0.2cm}
0 & b_1y_0 & 0 \\
 \vspace{0.2cm} 0 & 0 & 0\end{array}\right),
$$ with $V(t)=(u_1(t), u_2(t), u_3(t))^T$.

The characteristic equation of system (23) is given by
$\Delta(\lambda, \tau_1, q_2)=0$, where $\Delta(\lambda, \tau_1,
q_2)$ is function (16). We consider $\tau_1=\tau_{11}+\mu$,
$\mu\in\rr$ and $C^{1}=C^{1}([-\tau_{11},0],\mathbb{C}^{2})$. We
regard $\mu$ as the bifurcation parameter. Then, for $\phi\in
C^{1}$, we define

\begin{equation*}
\ca_1(\mu)\phi=\left\{\begin{array}{rl}
\frac{d\phi(\theta)}{d\theta},&\theta\in [-\tau_{11},0)\\
-\int\limits_{-\tau_{11}}^{0}d\eta(t,\mu)\phi(t),&
\theta=0\end{array} \right.
\end{equation*}and

\begin{equation*}
R_1(\mu)\phi=\left\{\begin{array}{rl} (0, 0, 0)^T, &\theta\in [-\tau_{11},0)\\
(-a_2\phi_1(0)\phi_2(0), b_{1}\phi_{3}(0)\phi_{2}(-\tau_{11}),
0)^T,& \theta=0
\end{array}\right.,
\end{equation*}
where

\begin{equation*}
\eta(\theta,\mu)=\left\{\begin{array}{rl}
A,&\theta=0\\
C_{1}\delta(\theta+\tau_{11}),& \theta\in [-\tau_{11},0).
\end{array} \right.
\end{equation*}\\ Then, we can rewrite (23) in
the following vector form

\begin{equation*}\dot{U}_{t}=\ca_1(\mu)U_{t}+R_1U_{t},
\end{equation*}where
\begin{equation*}
U_{t}=U(t+\theta),\quad \theta\in [-\tau_{11},0].
\end{equation*}
For $\psi \in C^{1}([0,\tau_{11}],\mathbb{C}^{2})$, the adjoint
operator $\ca_1^{\ast}$ of $\ca_1$ is defined as

\begin{equation*}\ca_1^{\ast}\psi(s)=\left\{\begin{array}{rl}
-\frac{d\psi(s)}{ds},& s\in (0,\tau_{11}] \\
\int\limits_{-\tau_{11}}^{0}d\eta^{T}(t,0)\psi(-t), & s=0.
\end{array}\right.
\end{equation*}

For $\phi\in C([-\tau_{11},0],\mathbb{C}^{2})$ and $\psi \in
C([0,\tau_{11}],\mathbb{C}^{2})$ we define the bilinear form

\begin{equation}<\phi,\psi>=\overline{\psi}^{T}(0)\phi(0)-
\int\limits_{-\tau_{11}}^{0}\int\limits_{\xi=0}^{\theta}
\overline{\psi}^{T}(\xi-\theta) d\eta(\theta)\phi(\xi)d \xi,
\end{equation}
where $\eta(\theta)=\eta(\theta,0)$.

\begin{proposition}\label{prop6}
\begin{itemize}
\item[(i)]The eigenvector  of $\ca_{1}(0)$ corresponding to
eigenvalue $\lambda_{1}=i\omega_{01}$ is given by
\begin{equation*}
h(\theta)=(v_{1},v_{2}, v_3)^{T}e^{\lambda_{1}\theta}, \quad
\theta\in[-\tau_{11},0],
\end{equation*} where

\begin{equation*}v_{1}=(\lambda_1+q_2)(\lambda_1+b_3-b_1y_0e^{\lambda_2\tau_{11}})-q_2b_1x_0,
v_{2}=-b_2(\lambda_{1}+q_2),  v_3=-b_2q_2
\end{equation*}
and $\lambda_{2}=\overline{\lambda_{1}}$; \item[(ii)]The
eigenvector  of $\ca_1^{\ast}$ corresponding to  eigenvalue
$\lambda_{2}$ is
\begin{equation*}h^{\ast}(s)=(w_{1},w_{2}, w_3)^{T}e^{\lambda_{1}s},
\quad s\in [0,\infty)
\end{equation*}
where
\begin{equation*}
\begin{split}
& w_{1}=\frac{f_1}{\eta}, w_{2}=\frac{1}{\eta},
w_{3}=\frac{f_3}{\eta},  f_1=-\frac{b_{2}}{\lambda_{2}},
f_3=\frac{b_{1}x_0}{\lambda_{2}+q_2}\\
& \eta=f_1\overline{v_{1}}+
\overline{v_{2}}(1-\frac{b_1y_0}{\lambda_2^2}(1-e^{\lambda_{1}\tau_{11}}-
\lambda_2\tau_{11}b_{2}e^{\lambda_{1}\tau_{11}}))+f_3\overline{v_{3}};
\end{split}
\end{equation*}
\item[(iii)]With respect to (24), we have:
\begin{equation*}<h^{\ast},h>=1,\quad
<h^{\ast},\overline{h}>=<\overline h^{\ast},h>=0, \quad <\overline
h^{\ast},\overline{h}>=1.
\end{equation*}
\end{itemize}
\end{proposition}

Using the approach of [2], [9] we next compute the coordinates of
the center manifold $\Omega_{0}$ at $\mu=0$. Let
$X_{t}=X(t+\theta), \theta\in [-\tau_{11},0)$ be the solution of
system (3) when $\mu =0$.\\
Define

\begin{equation*}
z(t)=<h^{\ast},X_{t}>\quad  w(t,\theta)=X_{t}- 2 Re
(z(t)h(\theta)).
\end{equation*}

On the center  manifold $\Omega_{0}$, we have
\begin{equation*}w(t,\theta)=w(z(t),\overline{z}(t),\theta)=
w_{20}(\theta)\frac{z^{2}}{2}+w_{11}(\theta)
z\overline{z}+w_{02}(\theta)\frac{\overline{z}^{2}}{2}+\cdots
\end{equation*}
where $z$ and $\overline{z}$ are the local coordinates of the
center manifold $\Omega_{0}$ in the direction of $h$ and
$h^{\ast}$, respectively.\\
For the solution $X_{t}\in\Omega_{0}$, we have:

\begin{equation*}\dot{z}(t)=\lambda_{1}z(t)+g(z(t),\overline{z(t)}),
\end{equation*}where
\begin{equation*}g(z,\overline{z})=g_{20}\frac{z^{2}}{2}+g_{11}z\overline{z}+g_{21}
\frac{\overline{z}^{2}}{2}+g_{21}\frac{z^{2}\overline{z}}{2}
\end{equation*}

\begin{proposition}For the system (23), the coefficients $g_{20}, g_{11},
 g_{02}, g_{21}$ and the functions $w_{20}(\theta), w_{11}(\theta), w_{02}(\theta)$
are given by
\begin{equation}
\begin{split}
& g_{20}=\overline w_1f_{120}+\overline w_2f_{220}+\overline
w_3f_{320},\quad
g_{11}=\overline w_1f_{111}+\overline w_2f_{211}+\overline w_3f_{311},\\
& g_{02}=\overline w_1f_{102}+\overline w_2f_{202}+\overline
w_3f_{302},\quad g_{21}=\overline w_1f_{121}+\overline
w_2f_{221}+\overline w_3f_{321}.
\end{split}
\end{equation}
where
\begin{equation*}
\begin{split}
&  f_{120}=-2a_2v_{1}v_{2}, f_{111}=-2a_2Re(v_1\overline{v_{2}}),
f_{102}=\overline
f_{120},\\
& f_{220}=2b_{1}v_2v_{3}e^{\lambda_{2}\tau_{11}},
f_{211}=2b_{1}Re(\overline v_2v_{3}e^{\lambda_{2}\tau_{11}}),
f_{202}=\overline f_{220}\\
&f_{320}=f_{311}=f_{302}=0
\end{split}
\end{equation*}

\begin{equation*}
\begin{split}
& f_{121}=-a_{2}(2v_{1}w_{211}(0)+\overline{v_{1}}w_{220}(0)+2v_{2}w_{111}(0)+\overline{v_{2}}w_{120}(0))\\
&
f_{221}=b_{1}(2v_{3}w_{211}(-\tau_{11})+\overline{v_{3}}w_{220}(-\tau_{11})+2v_2e^{\lambda_{2}\tau_{11}}w_{311}(0)+
\overline{v_{2}}e^{\lambda_{1}\tau_{11}}w_{320}(0))\\
& f_{321}=0;
\end{split}
\end{equation*}
and
\begin{equation*}
w_{20}(\theta)=(w_{120}(\theta), w_{220}(\theta),
w_{320}(\theta))^T, \quad w_{11}(\theta)=(w_{111}(\theta),
w_{211}(\theta), w_{311}(\theta))^T,
\end{equation*}

\begin{equation*}
\begin{split}
&w_{20}(\theta)=-\frac{g_{20}}{\lambda_{1}}h(0)e^{\lambda_{1}\theta}
-\frac{\overline{g_{20}}}{3\lambda_{1}}\overline{h(0)}e^{\lambda_{2}\theta}+
E_{1}e^{2\lambda_{1}\theta}\\
&w_{11}(\theta)=\frac{g_{11}}{\lambda_{1}}h(0)e^{\lambda_{1}\theta}-
\frac{\overline{g_{11}}}{\lambda_{1}}\overline{h(0)}e^{\lambda_{2}\theta}+E_{2},
\end{split}
\end{equation*}
and
\begin{equation*}
E_{1}=(E_{11}, E_{12}, E_{13})^{T}, \quad E_{2}=(E_{21}, E_{22},
E_{23})^{T},
\end{equation*}
where
\begin{equation*}
\begin{split}
&E_{11}=\frac{a_2x_0E_{12}+f_{120}}{2\lambda_{1}},
E_{12}=\frac{2\lambda_1f_{220}-b_2f_{120}}{2\lambda_1(2\lambda_1+b_3-b_1y_0e^{2\lambda_1\tau_{11}})-a_2b_2x_0},\\
&E_{13}=\frac{q_2E_{11}}{2\lambda_1+q_2},
E_{21}=-\frac{(b_3-b_1y_0)E_{22}-f_{211}}{b_2},
E_{22}=\frac{f_{111}}{a_2x_0}, E_{23}=E_{21}.
\end{split}
\end{equation*}
\end{proposition}

We can explicitly compute the following quantities $C_{11}(0)$,
$\mu_{21}$, $T_{21}$, $\beta_{21}$:
\begin{equation}
\begin{split}&C_{11}(0)=\frac{i}{2\omega_{11}}(g_{20}g_{11}-2|g_{11}|^{2}-\frac{1}{3}
|g_{02}|^{2})+\frac{g_{21}}{2}\\
&\mu_{21}=-\frac{Re\, C_{11}(0)}{Re \,\lambda'(0)}\\
&T_{21}=-\frac{Im\,C_{11}(0)+\mu_{21}Im\,\lambda'(0)
}{\omega_{11}}\\
&\beta_{21}=2Re\, (C_{11}(0)).
\end{split}
\end{equation}

In summary, this leads to the following result:
\begin{theorem}In formulas (26), $\mu_{21}$ determines the directions of the Hopf
bifurcations: if $\mu_{21}>0(<0)$ the Hopf bifurcation is
supercritical (subcritical) and the bifurcating periodic solutions
exist for $\tau_{1}>\tau_{11}(<\tau_{11})$; $\beta_{21}$
determines the stability of the bifurcation periodic solutions:
the solutions are orbitally stable (unstable) if $\beta_{21}< 0 (>
0)$ and $T_{21}$ determines the periodic solutions: the period
increases (decreases)
 if $T_{21}>0 (<0)$. \\
In $(26)$ $Re (\lambda'(0))$ and $Im(\lambda'(0))$ are given by
\begin{equation*}
\begin{split}
&Re(\lambda'(0))=Re(\frac{d\lambda}{d\tau_{1}})_{\lambda=i\omega_{01},
\tau=\tau_{11}}\\
&
Im(\lambda'(0))=Im(\frac{d\lambda}{d\tau_{1}})_{\lambda=i\omega_{01},
\tau=\tau_{11}}
\end{split}
\end{equation*}
where
\begin{equation*}
\frac{d\lambda}{d\tau_{1}}=\frac{(r_1\lambda^2+r_0\lambda-r_1)
e^{-\lambda\tau_{1}}}{3\lambda^2+2p_2\lambda+p_1-(r_1\lambda+r_0)\tau_1}.
\end{equation*}
\end{theorem}

\medskip
\medskip

\section*{\normalsize\bf 4. Numerical simulations.}

\hspace{0.6cm}For the numerical simulations we use Maple 9.5. In
this section, we consider system (6) with $a_1=2.5, a_2=1$,
$b_1=1$, $b_2=0.4$, $b_3=0.95$, $b_{4}=2$. We obtain: $x_{0}=
 0.1524390244$, $y_{0}=2.5$.

In the first case, $k_1(s)=\delta(s-\tau_1)$,
$k_2(s)=\delta(s-\tau_2)$, for $\tau_2=0.01$, we have: $\omega_{0}
= 0.6124295863$, $\mu_2\!=\! 630.5712553$, $\beta_2\!=\!
125.5070607$, $T_2\!=\! 10.25944116$, $\tau_{10}\!=9.541873607\!$.
Then the Hopf bifurcation is supercritical and the bifurcating
periodic solutions exist for $\tau>\tau_{10}$; the solutions are
orbitally unstable and the period of the solution increases. The
waveforms are displayed in Fig1 and Fig2 and the phase plane
diagram of the state variables $x(t)$, $y(t)$ is displayed in
Fig3:

\begin{center}\begin{tabular}{ccc}
%\hline
\includegraphics[width=6cm]{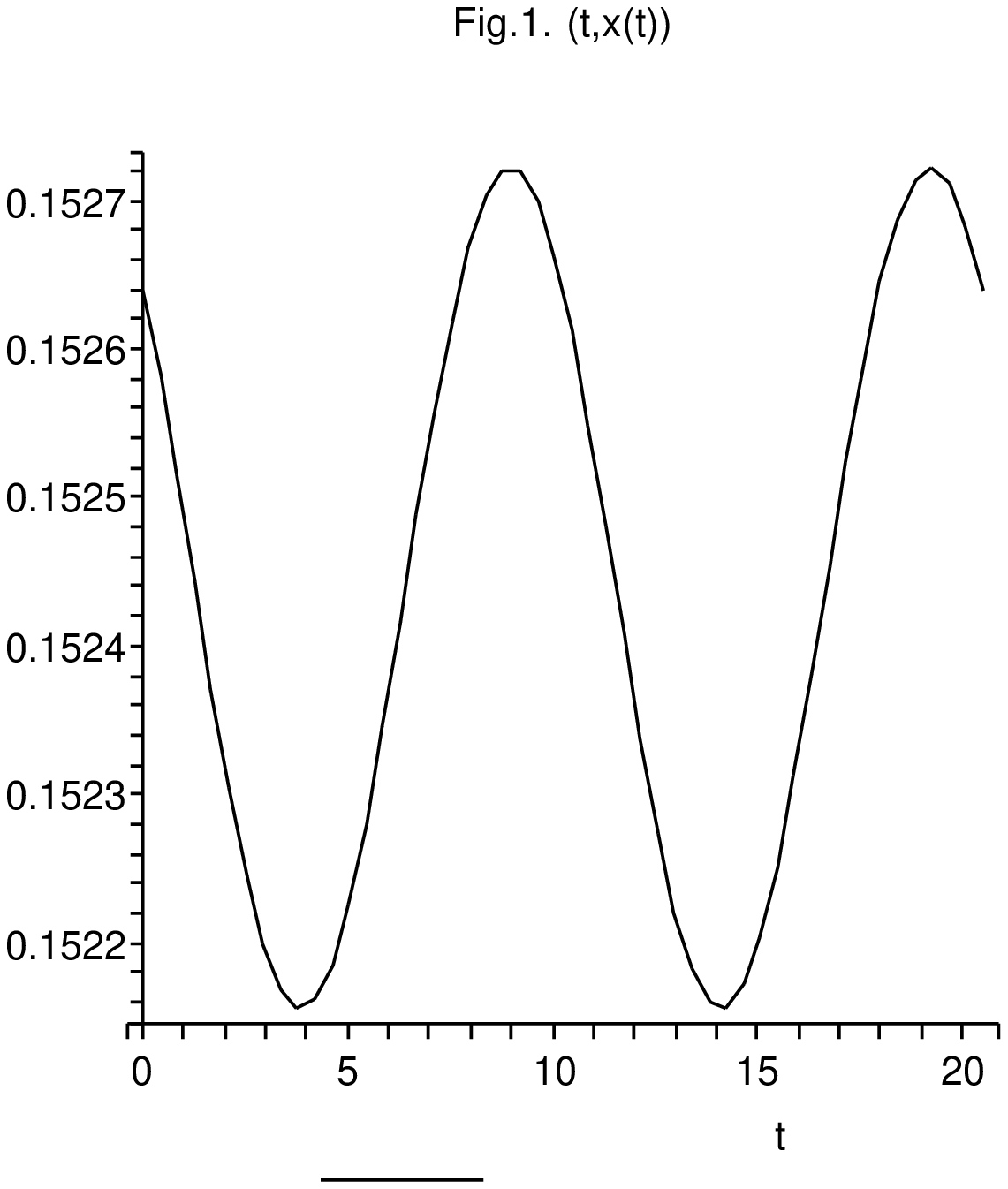}
\includegraphics[width=6cm]{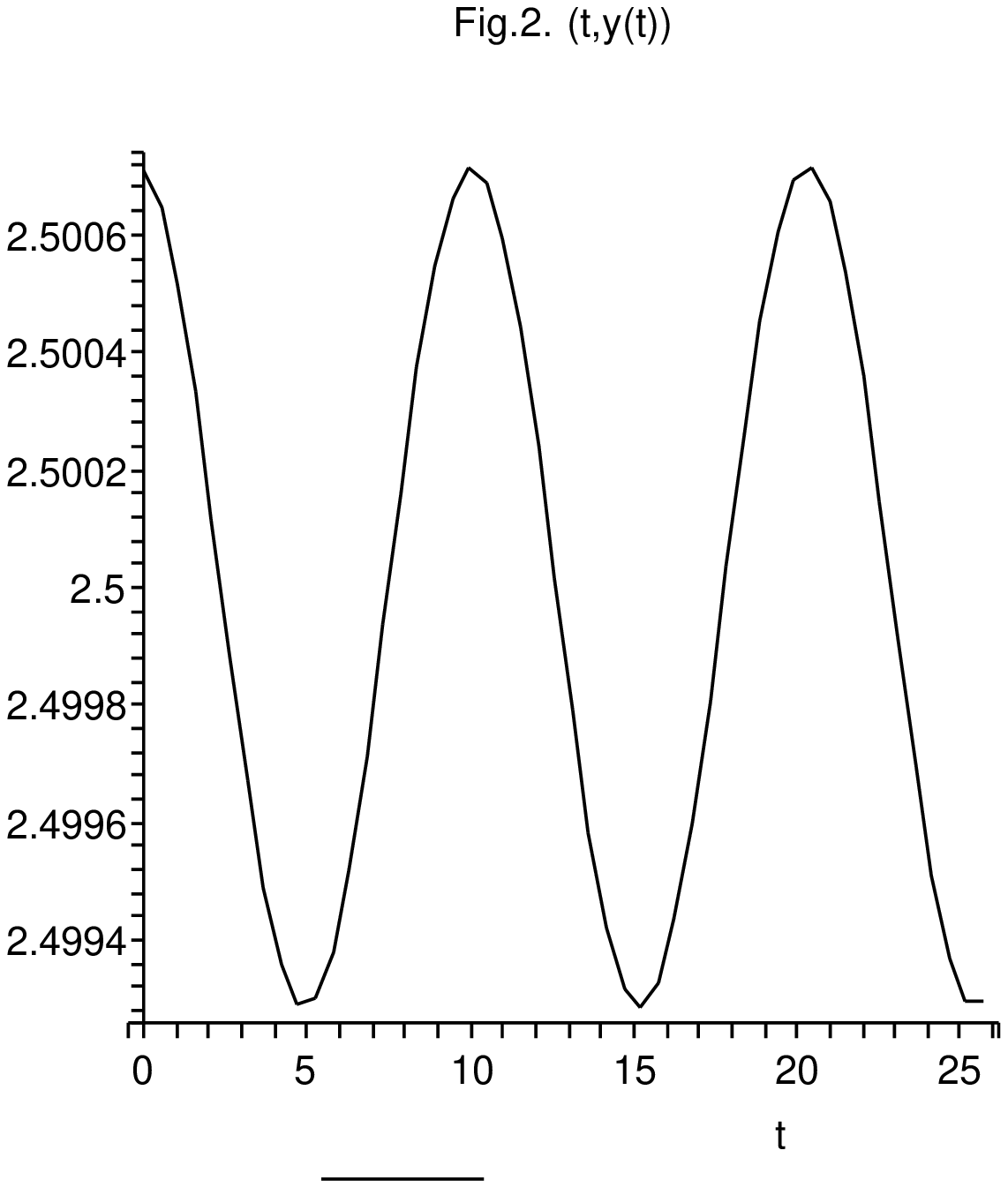}

\end{tabular}
\end{center}
\begin{center}\begin{tabular}{cc}
%\hline

\includegraphics[width=6cm]{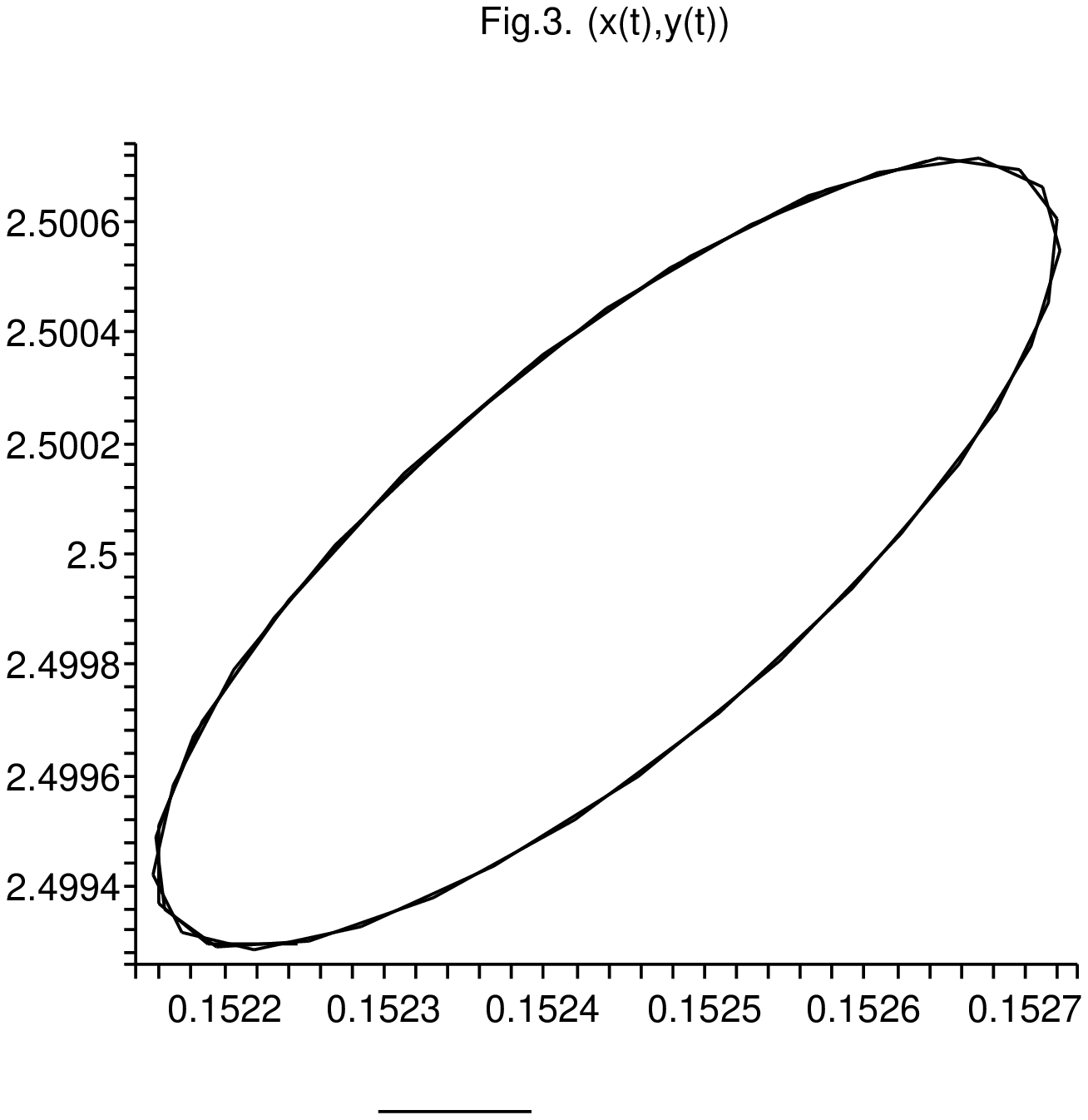}

\end{tabular}
\end{center}

\medskip

In the second case, $k_1(s)=\delta(s-\tau_1)$,
$k_2(s)=q_2e^{-q_2s}$ for $q_2=0.1$, we have:  $\omega_{01} =
0.2235621332$, $\mu_{21}\!=\!7.926079992$, $\beta_{21}\!=\!
0.04097046568$, $T_{21}\!=\! 0.3275619874$,
$\tau_{11}^{\ast}\!=10.38589492$. Then the Hopf bifurcation is
supercritical and the bifurcating periodic solutions exist for
$\tau_1>\tau_{11}^{\ast}$; the solutions are orbitally unstable
and the period of the solution increases. The waveforms are
displayed in Fig5 and Fig6 and the phase plane diagram of the
state variables $x(t)$, $y(t)$ is displayed in Fig7:
\begin{center}\begin{tabular}{ccc}
%\hline
\includegraphics[width=6cm]{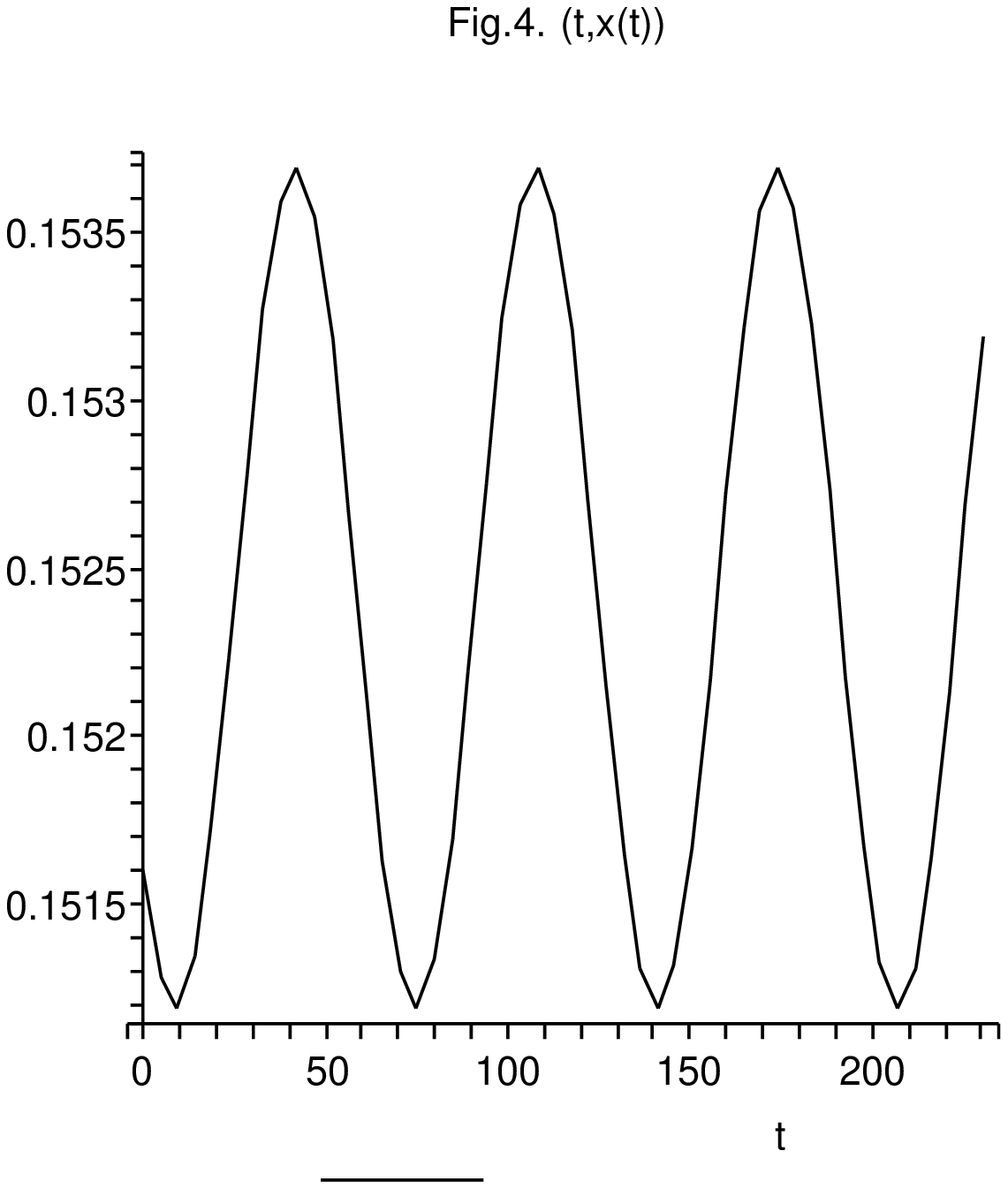}
\includegraphics[width=6cm]{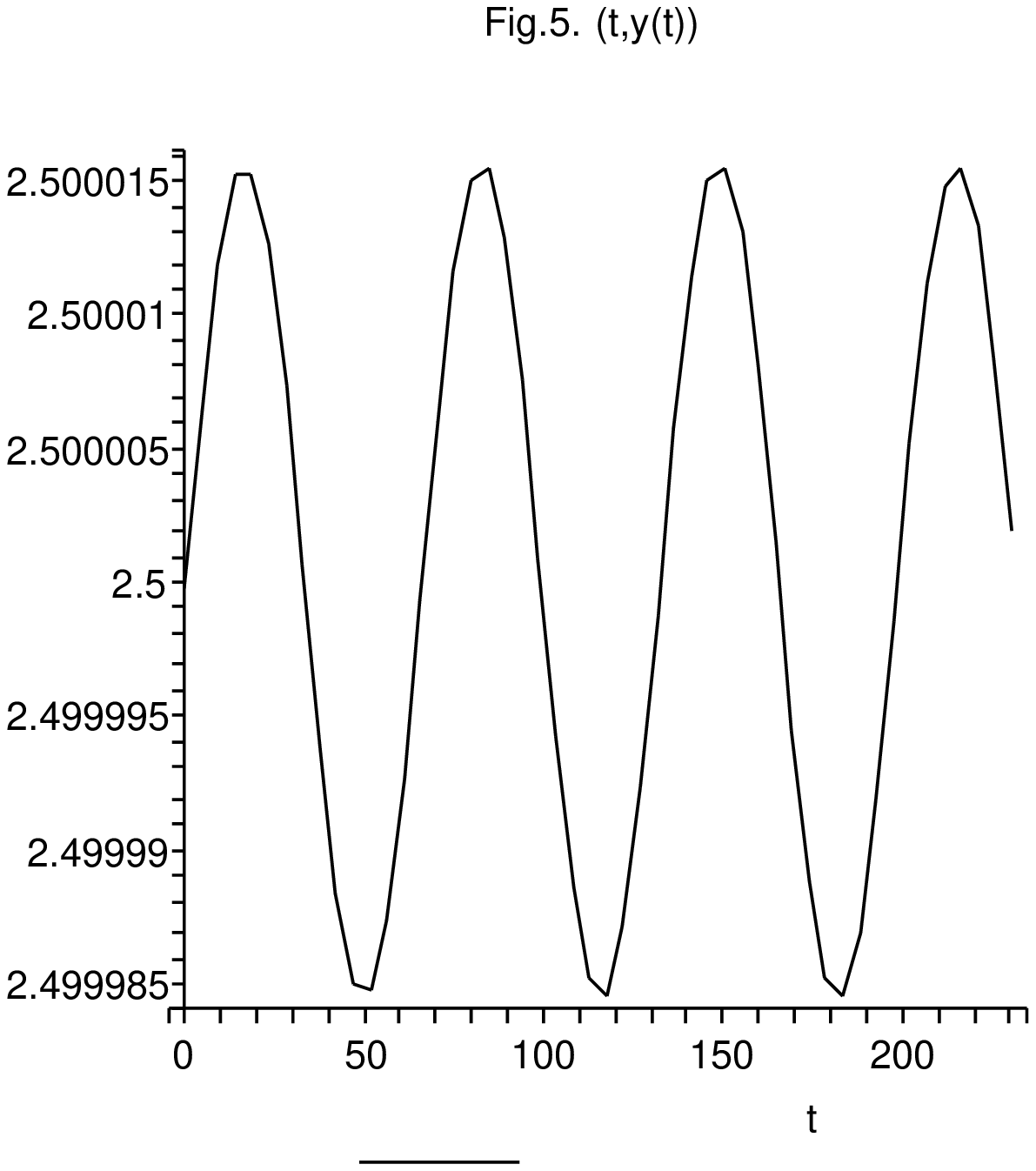}

\end{tabular}
\end{center}
\begin{center}\begin{tabular}{cc}
%\hline

\includegraphics[width=6cm]{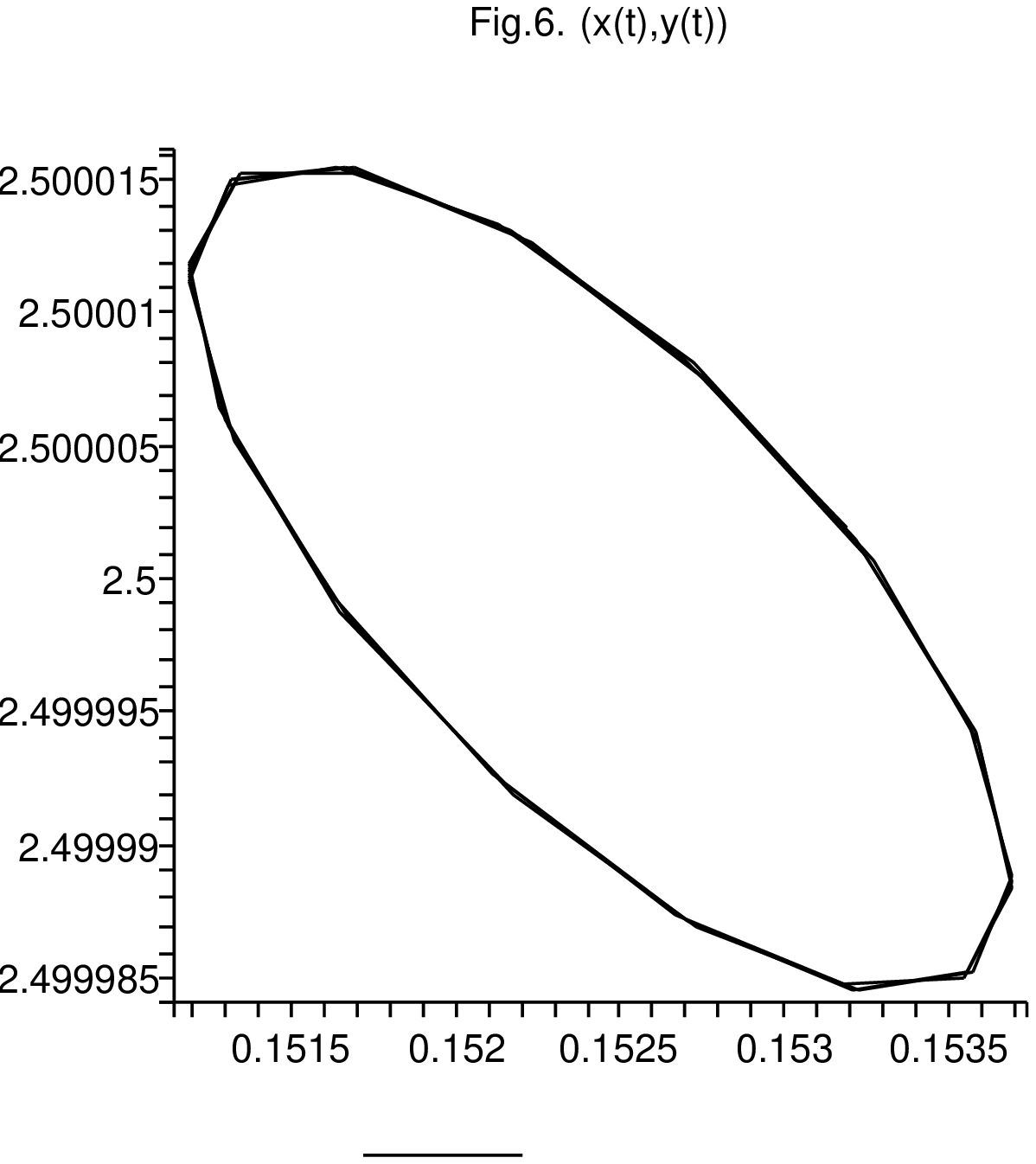}

\end{tabular}
\end{center}

\medskip

For $q_2=0.1$, we have:  $\omega_{01} = 0.9506753825$,
$\mu_{21}\!=\!-0.6058263333$, $\beta_{21}\!=\! -0.001118156944$,
$T_{21}\!=\! -0.07864963978$, $\tau_{11}\!=23.03933807$. Then the
Hopf bifurcation is subcritical and the bifurcating periodic
solutions exist for $\tau_1>\tau_{11}$; the solutions are
orbitally stable and the period of the solution decreases. The
waveforms are displayed in Fig5 and Fig6 and the phase plane
diagram of the state variables $x(t)$, $y(t)$ is displayed in Fig7
and Fig8:

\begin{center}\begin{tabular}{ccc}
%\hline
\includegraphics[width=6cm]{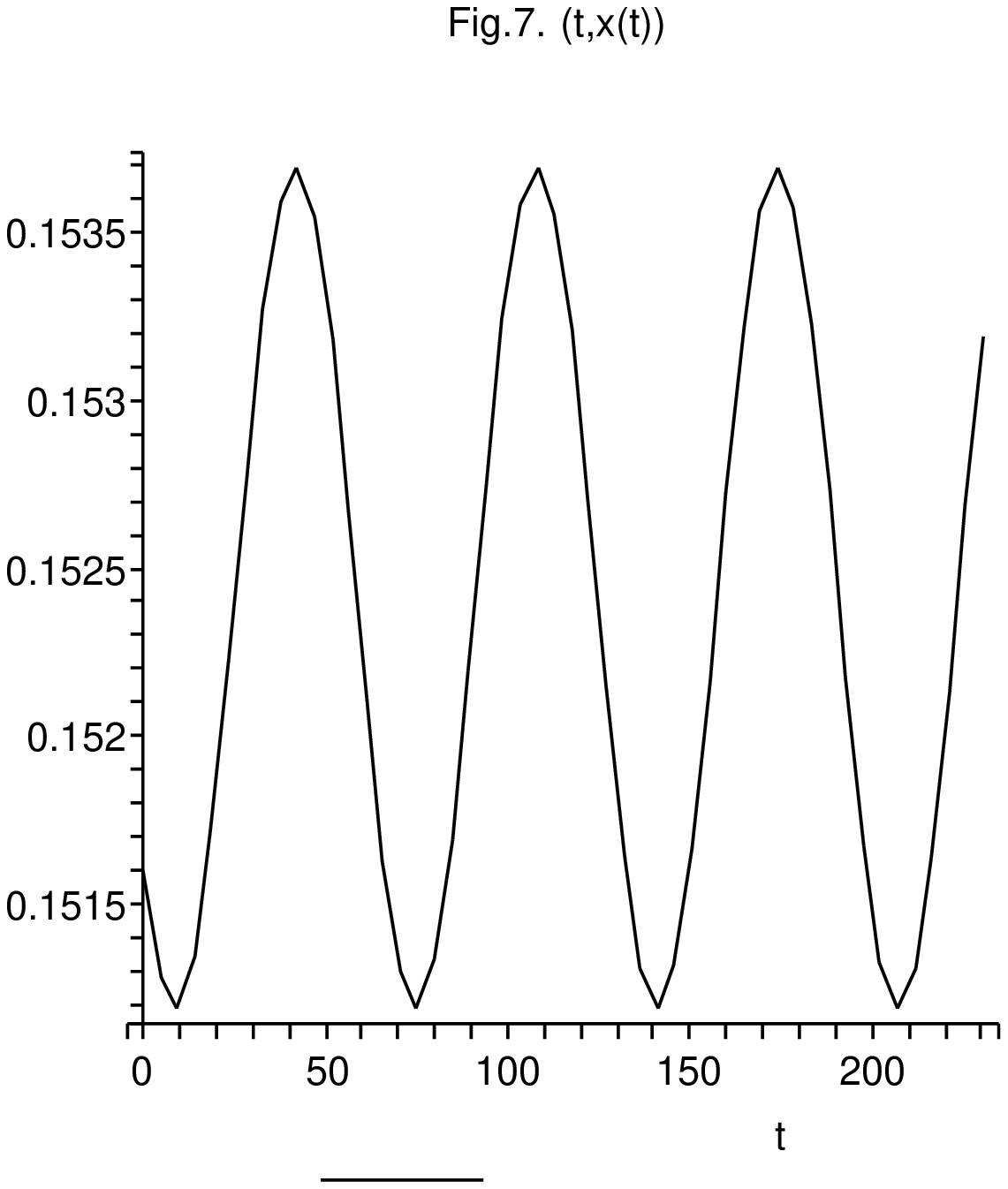}
\includegraphics[width=6cm]{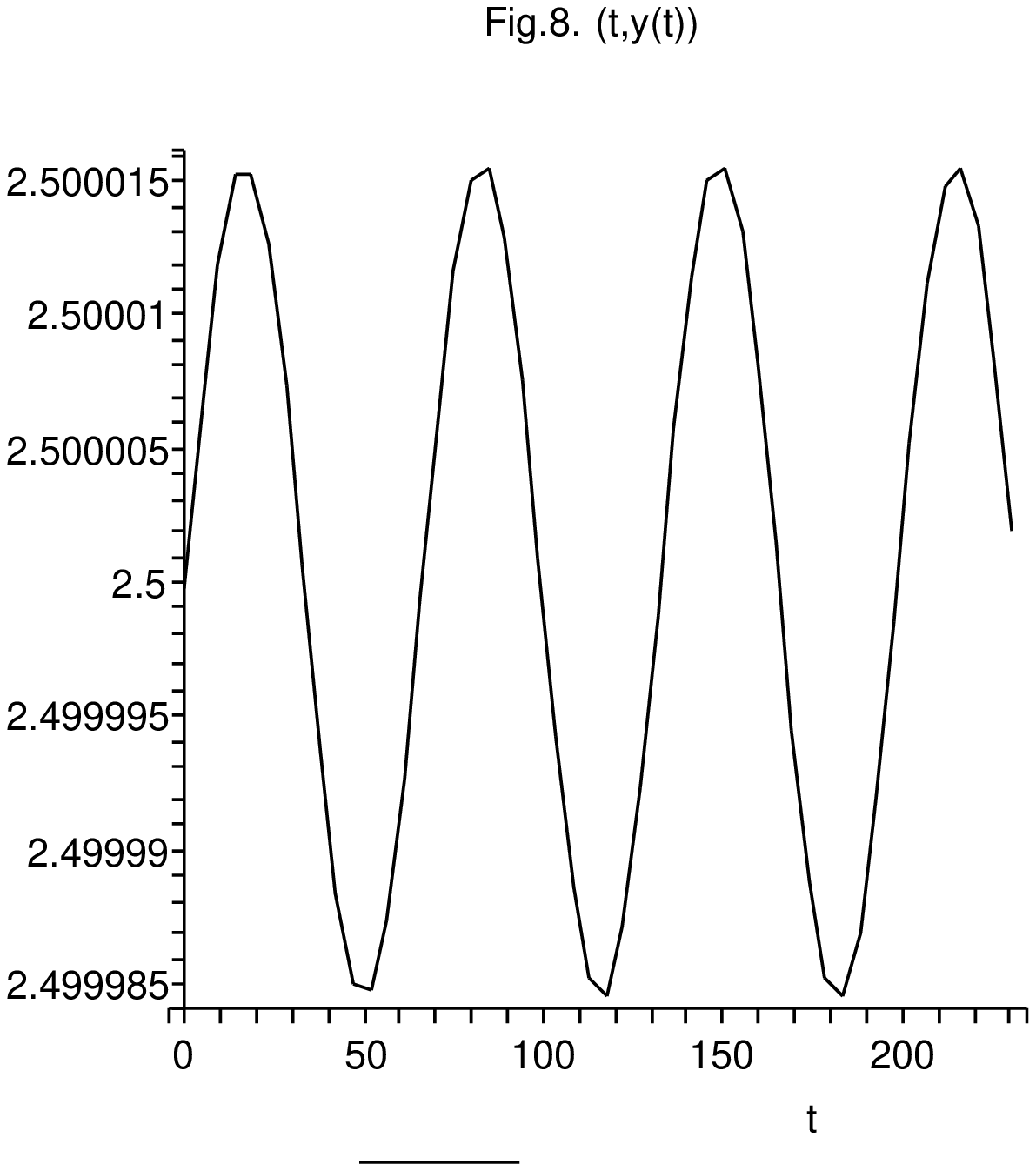}

\end{tabular}
\end{center}
\begin{center}\begin{tabular}{cc}
%\hline

\includegraphics[width=6cm]{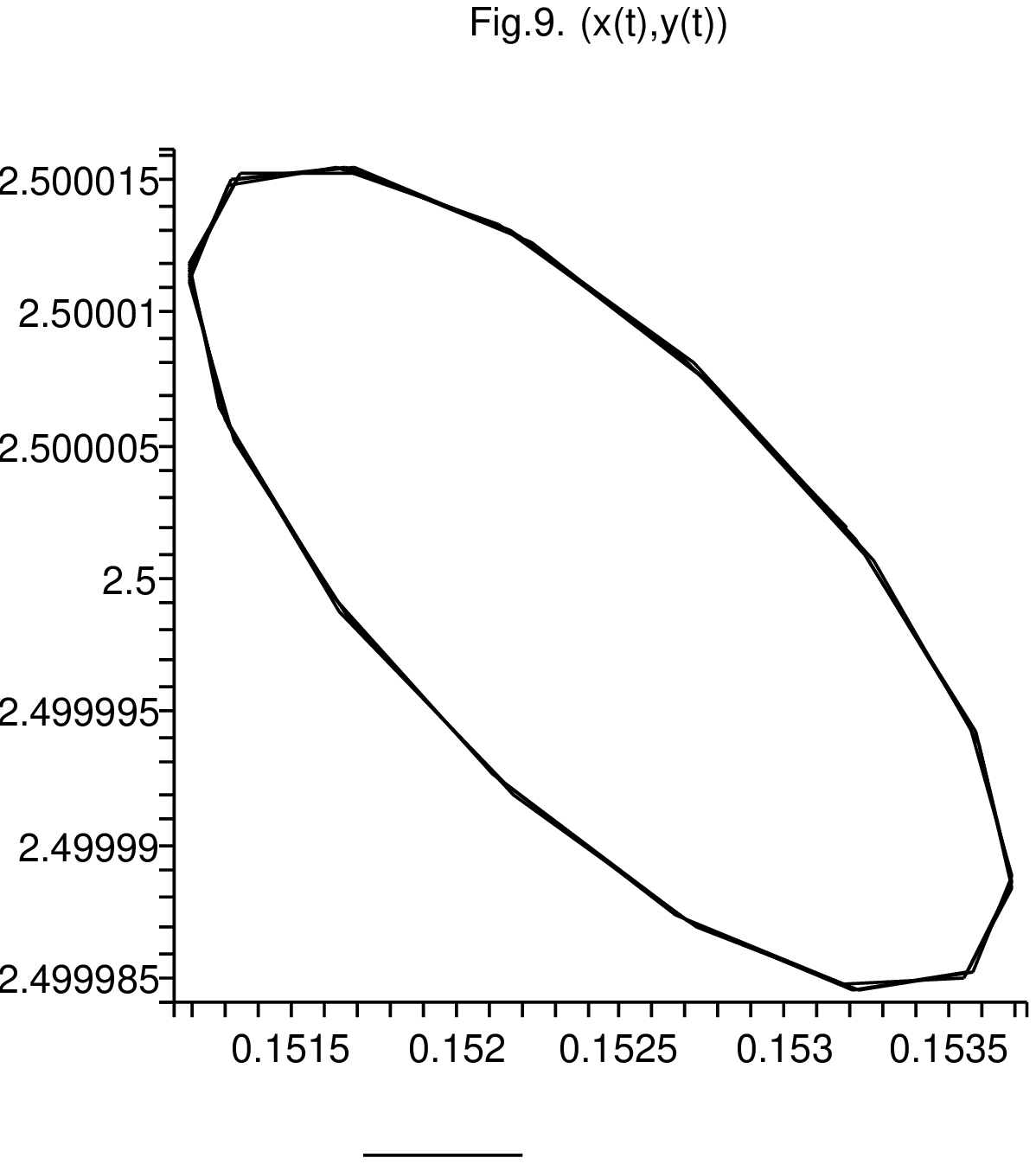}

\end{tabular}
\end{center}

\section*{\normalsize\bf 5. Conclusions.}

This paper was focused on mathematical analysis of a model which
describes the interaction between immune system and the tumor
cells. The model is an improved one by using the delay kernel.
Taking the average time delay as a parameter, it has been proved
that the Hopf bifurcation occurs when this parameter passes
through a critical value. In a future work it will be studied the
mathematical aspects of the effect of immunotherapy on the
development of the cancer.

\end{document}